\documentclass[reqno]{amsart}
\pdfoutput=1
\usepackage{amssymb, enumerate, mathrsfs, bbm, verbatim, array}
\usepackage[centredisplay, PostScript=dvips]{diagrams}
\usepackage[pdftex]{graphicx}
\usepackage[initials]{amsrefs}

\title{The Automorphism Group of a Metropolis-Rota Implication Algebra}
\date{2009, February 6}
\author{Colin G. Bailey}
\address{School of Mathematics,  Statistics \& Operations Research\\
Victoria University of Wellington\\
Wellington, New Zealand\\
}
\email{Colin.Bailey@vuw.ac.nz}
\author{Joseph S. Oliveira}
\address{
Pacific Northwest National Laboratories\\
Richland, Washington\\
U.S.A.}
\email{Joseph.Oliveira@pnl.gov}
\let\rsf\mathscr

\def\caret{\mathbin{\hat{\hphantom{n}}}} 

\def\one{{\mathbf 1}}
\def\aut{\operatorname{Aut}}
\def\incl{\operatorname{incl}}
\def\inn{\operatorname{Inn}}
\def\leftGen{{[\kern-1.1pt[}}
\def\rightGen{{]\kern-1.1pt]}}

\providecommand{\meet}{\mathbin{\wedge}}
\providecommand{\join}{\mathbin{\vee}}
\newcommand{\comp}[1]{\overline{#1}}

\ifx\upharpoonright\undefined
     \def\restrict{\hbox{\rm\kern0.166em\accent"12\kern-0.536em$\vert$\kern0.3em}}%
  \else
     \def\restrict{\upharpoonright}%
  \fi
\makeatletter
\def\twoSet#1#2{\left\{%
\vphantom{#2}#1\thinspace\right|\nolinebreak[3]\left.%
  #2%
  \vphantom{#1}%
  \right\}%
}
\def\oneSet#1{\left\lbrace#1\right\rbrace}

\newif\if@nstr
\def\setstrfalse{\let\if@nstr=\iffalse}
\def\setstrtrue{\let\if@nstr=\iftrue}
\def\@nstr #1#2{
\def\@@nstr ##1#1##2##3\@@nstr{\ifx
\@nstr ##2\setstrfalse \else \setstrtrue \fi }
\@@nstr #2#1\@nstr \@@nstr}
\def\@separate#1|#2@{\setFront{#1}\setBack{#2}}
\def\lb#1\rb{\@nstr|{#1} \if@nstr \@separate#1 @ \twoSet{\@setFront}{\@setBack}%
\else \@separate |{#1 }@ \oneSet{\@setBack}\fi%
}
\def\setFront#1{\def\@setFront{#1}}
\def\setBack#1{\def\@setBack{#1}}
\def\Set#1{\lb{#1}\rb}
\def\oneBrk#1{\left\langle#1\right\rangle}
\def\twoBrk#1#2{\left\langle%
\vphantom{#2}#1\thinspace\right|\nolinebreak[3]\left.%
  #2%
  \vphantom{#1}%
  \right\rangle%
}
\def\brk<#1>{\@nstr|{#1} \if@nstr \@separate#1 @ \twoBrk{\@setFront}{\@setBack}%
\else \@separate |{#1 }@ \oneBrk{\@setBack}\fi%
}
\def\thmref#1{\normalfont{theorem}~\ref{#1}}
\def\lemref#1{\normalfont{lemma}~\ref{#1}}

\def\corref#1{\normalfont{corollary}~\ref{#1}}

\theoremstyle{plain}
\newtheorem{thm}{Theorem}[section]
\newtheorem{lem}[thm]{Lemma}
\newtheorem{cor}[thm]{Corollary}
\newtheorem{prop}[thm]{Proposition}
\newtheorem{defn}[thm]{Definition}

\theoremstyle{remark}
{}
{}
{}
{}

\@ifpackageloaded{bbm}{\newcommand{\Z}{{\mathbbm{Z}}}

}{
\newcommand{\Z}{{\mathbb{Z}}}

}
\makeatother

\begin{document}

\begin{abstract}
	We discuss the group of automorphisms of a general MR-algebra. We 
	develop several functors between implication algebras and cubic 
	algebras. These allow us to generalize the notion of inner 
	automorphism. We then show that this group is always isomorphic to the 
	group   of inner automorphisms of a filter algebra. 
\end{abstract}
\maketitle

\section{Introduction}
Cubic implication algebras are an algebraic generalization of the algebra of faces of 
an $n$-cube as introduced by Metropolis \& Rota in \cite{MR:cubes}. 
From their paper and further work of the authors (\cites{BO:eq, 
BO:fil}) much is known about the structure of cubic implication algebras. A 
survey of this material may be found in \cite{BO:surv}.

The group of automorphisms of a face poset of 
an $n$-cube is well-known to be $\Z_{2}^{n}\rtimes S_{n}$. 
We will show that in every cubic implication algebra there is a subgroup of definable automorphisms, 
known as  \emph{inner automorphisms}, that corresponds to the $\Z_{2}^{n}$ portion 
of this group.  

Associated with any cubic implication algebra is a minimal enveloping 
Metropolis-Rota implication algebra (usually abbreviated as MR-algebra). 
Notions such as congruences and automorphisms 
lift from the cubic implication algebra to its envelope. 
Thus a major part of characterizing the automorphism group of a cubic 
implication
algebra is the special case of MR-algebras. Earlier work 
(\cite{BO:fil}) has dealt 
with proper subclasses -- interval algebras and filter algebras. 
Herein we are interested in automorphism groups of arbitrary MR-algebras. 
By careful construction of subalgebras we are able to lift results 
from filter and interval algebras to the general case. 

We begin with some definitions and basic results from \cite{BO:eq}. 
\begin{defn}
	A \emph{cubic implication algebra} is a join semi-lattice with one and a binary 
	operation $\Delta$ satisfying the following axioms:
	\begin{enumerate}[a.]
		\item  if $x\le y$ then $\Delta(y, x)\join x = y$;
		
		\item  if $x\le y\le z$ then $\Delta(z, \Delta(y, x))=\Delta(\Delta(z, 
		y), \Delta(z, x))$;
		
		\item  if $x\le y$ then $\Delta(y, \Delta(y, x))=x$;
		
		\item  if $x\le y\le z$ then $\Delta(z, x)\le \Delta(z, y)$;
		
		\item[] Let $x\to y=\Delta(\one, \Delta(x\join y, y))\join y$ for any $x$, $y$ 
		in $\mathcal L$. Then:
		
		\item  $(x\to y)\to y=x\join y$;
		
		\item  $x\to (y\to z)=y\to (x\to z)$;
	\end{enumerate}
\end{defn} 

\begin{defn}
	An \emph{MR-algebra} is a cubic implication algebra satisfying the MR-axiom:\\
	if $a, b<x$ then 
	\begin{gather*}
		\Delta(x, a)\join b<x\text{ iff }a\meet b\text{ does not exist.}
	\end{gather*}
\end{defn}

\begin{defn}\label{def:caret}
    Let $\mathcal L$ be a cubic implication algebra. Then for any $x, y\in\mathcal L$ 
    we define the (partial) operations $\caret$ (\emph{caret}) 
    and $*$ by:
	\begin{enumerate}[(a)]
		\item
    $$
        x\caret y=x\meet\Delta(x\join y, y)
    $$
    whenever this meet exists. 
	\item 
	$$
        x * y=x\join\Delta(x\join y, y)
    $$
	\end{enumerate}
\end{defn}

\begin{lem}
    Every cubic implication algebra is an implication algebra.
\end{lem}
\begin{proof}
    It suffices to note that $x\to x= \Delta(\one, \Delta(x, x))\join 
    x= \Delta(\one, x)\join x= \one$.
\end{proof}

\begin{lem}
    If $\mathcal L$ is a cubic implication algebra then 
    $\mathcal L$ is an MR-algebra iff the caret operation is total.
\end{lem}
\begin{proof}
    If $\mathcal L$ is an MR-algebra, then for any $a,b$ we have 
    $a\join b=a\join b$ and so $a\meet\Delta(a\join b,b)= a\caret b$ 
    exists. 
    
    Conversely, suppose caret is total. If $a\meet b$ exists, then
    $x\geq a\join\Delta(x, b)\geq(a\meet b)\join\Delta(x, a\meet b)=x$.
    
    Now suppose that $a\join\Delta(x, b)=x$. There are two cases 
    \begin{enumerate}[-]
	    \item  if $a\join b$ is one of $a$ or $b$,  then $a$ and $b$ are 
	    comparable and the meet clearly exists.
    
	    \item  Otherwise $a, b<a\join b$. By \cite{BO:eq} 
	    theorem 4.3 we must have 
	    \begin{equation}\label{eq:one}
		    a\join\Delta(a\join b, b)=a\join b.
	    \end{equation}
	    Then we have
	    \begin{align*}			
		    a\caret\Delta(a\join b, b)&=a\meet\Delta(a\join\Delta(a\join b, b), 
		    \Delta(a\join b, b))&&\text{ by definition}\\
		    &=a\meet\Delta(a\join b, \Delta(a\join b, b))&&\text{ by 
		    \eqref{eq:one}}\\
		    &=a\meet b.
	    \end{align*}
    \end{enumerate}
\end{proof}

\begin{defn}
	Let $\mathcal L$ be a cubic implication algebra and $a, b\in\mathcal L$. Then
	\begin{align*}
		a\preccurlyeq b &\text{ iff }\Delta(a\join b, a)\le b\\
		a\sim b &\text{ iff }\Delta(a\join b, a)=b.
	\end{align*}
\end{defn}

\begin{lem}
	Let $\mathcal L$, $a$, $b$ be as in the definition above. Then
	$$
	a\preccurlyeq b\text{ iff }b=(b\join a)\meet(b\join\Delta(\one, a)).
	$$
\end{lem}
\begin{proof}
	See \cite{BO:eq} lemmas 2.7 and 2.12.
\end{proof}

\begin{prop}\label{prop:triv}
    Let $\mathcal L$ be a cubic implication algebra, and $p, q$ in $\mathcal L$ 
    are such that
    $p\preccurlyeq q$ and $p\meet q$ exists. Then $p\le q$.
\end{prop}
\begin{proof}
    We have $\Delta(p\join q, p)\le q$ as $p\preccurlyeq q$. Let
    $a=p\meet q$. Then
    \begin{align*}
        a&\le q\\
        \Delta(p\join q, a)&\le\Delta(p\join q, p)\le q\\
        \text{Hence }\qquad
        p\join q&=a\join\Delta(p\join q, a)\\
        &\le q
        \end{align*}
        and so $p\le q$.
\end{proof}

\begin{cor}\label{cor:triv}
    Let $\mathcal L$ be a cubic implication algebra, and $p, q$ in $\mathcal L$ 
    are such that
    $p\simeq q$ and $p\meet q$ exists. Then $p= q$.    
\end{cor}

Also from \cite{BO:eq} (lemma 2.7c for transitivity) we know that $\sim$ 
is an equivalence relation on all cubic implication algebras, and is a congruence 
on the structure $\brk<\mathcal L, \caret, *, \one>$. The quotient is 
naturally an implication algebra. 

As part of the representation theory in \cite{BO:eq} we had the following definitions 
and lemma:
\begin{defn}
	Let $\mathcal L$ be a cubic implication algebra. Then
	\begin{enumerate}[(a)]
		\item  $\mathcal L_{a}=\Set{\Delta(y, x) | a\le x\le y}$ for any 
		$a\in\mathcal L$. $\mathcal L_{a}$ is the \emph{localization} of $\mathcal 
		L$ at $a$. 
	
		\item  $k_{a}(y)=(\Delta(\one, y)\join a)a$ for any $y\in\mathcal L_{a}$.
	
		\item  $\ell_{a}(y)=y\join a$ for any $y\in\mathcal L_{a}$.
	\end{enumerate}
\end{defn}

\begin{lem}\label{lem:kl}\label{lem:onto}
	Let $\mathcal L$ be any cubic implication algebra. Then
	\begin{enumerate}[(a)]
		\item  $\mathcal L_{a}$ is an atomic MR-algebra, and hence 
		isomorphic to an interval algebra.
	
		\item  $k_{a}(x)\le\ell_{a}(x)$ for any $x\in\mathcal L_{a}$. 
	
		\item  If $x, y\in\mathcal L_{a}$ then 
			\begin{align*}
				x=y&\iff k_{a}(x)=k_{a}(y) \text{ and }\ell_{a}(x)=\ell_{a}(y)\\
				&\iff x\join a=y\join a\text{ and }x\join\Delta(\one, 
				a)=y\join\Delta(\one, a).
			\end{align*}		
		\item  If $p\geq q\geq a$ then there exists a unique $z\in\mathcal L_{a}$ 
		such that $\ell_{a}(z)=p$ and $k_{a}(z)=q$.
		
		\item $x\in\mathcal L_{a}$ iff $a\preccurlyeq x$. 
	\end{enumerate}
\end{lem}
\begin{proof}
	See \cite{BO:eq}.
\end{proof}

\section{Construction of cubic implication algebras}
There is a general construction of cubic implication algebras from implication 
algebras. 

Let $\mathcal I$ be an implication algebra. We define
$$
\rsf I(\mathcal I)=\Set{\brk<a, b> | a, b\in\mathcal I, a\join
b=1\text{ and }a\meet b\text{ exists}}
$$
ordered by 
$$
\brk<a, b>\le\brk<c, d>\text{ iff }a\le c\text{ and }b\le d. 
$$
This is a partial order that is an upper semi-lattice with join 
defined by
$$
\brk<a, b>\join\brk<c, d>=\brk<a\join c, b\join d>
$$
and a maximum element $\one=\brk<1, 1>$. 

We can also define a $\Delta$ function by
$$
\text{if }\brk<c, d>\le\brk<a, b>\text{ then }
\Delta(\brk<a, b>, \brk<c, d>)=\brk<a\meet(b\to d), b\meet(a\to c)>. 
$$

We note the natural embedding of $\mathcal I$ into $\rsf I(\mathcal 
I)$ given by
$$
e_{\mathcal I}(a)=\brk<1, a>. 
$$

Note that in an implication algebra $a\join b=\one$ iff $a\to 
b=b$ iff $b\to a=a$. 

Also $\Delta(\one, \bullet)$ is very simply defined as it is 
exactly $\brk<a, b>\mapsto\brk<b, a>$. 

It is not hard to show that $\rsf I(\mathcal I)$ is a cubic implication algebra, 
and is an MR-algebra iff $\mathcal I$ is a lattice. 

$\mathcal I$ is a lattice produces two cases -- either there is a 
least element or not. 

In the first of these cases $\mathcal I$ is a Boolean algebra and
$\rsf I(B)$ is naturally isomorphic to the algebra of closed intervals 
of $B$ -- we call these algebras \emph{interval algebras}. 

In the second case $\mathcal I$ is isomorphic to an ultrafilter of some 
Boolean algebra $B$ and we get a \emph{filter algebra} which can be 
embedded as an upwards closed MR-subalgebra of $\rsf I(B)$. 

Every interval algebra is also a filter algebra, as any Boolean 
algebra $B$ is an ultrafilter in $B\times\mathbf{2}$.  

MR-algebras that are isomorphic to filter algebras have an 
automorphism group that splits as a twisted product of a group of 
inner automorphisms with the group of automorphisms of the filter. The 
former group is isomorphic to a Boolean algebra that is  naturally definable 
from the filter. 
It is important in what follows to know that there are many MR-algebras 
that are filter algebras. In particular ones that are countably 
presented. 

\begin{defn}\label{def:presented}
	Let $\mathcal L$ be a cubic implication algebra. 
	\begin{enumerate}[(a)]
		\item Let $A\subseteq\mathcal L$. $A$ 
	is a \emph{presentation} of $\mathcal L$ iff 
	$$
		\mathcal L=\bigcup_{a\in A}\mathcal L_{a}. 
	$$
	
		\item $\mathcal L$ is \emph{countably presented} iff there is a 
		countable set presenting $\mathcal L$. 
	\end{enumerate}	
\end{defn}

\begin{thm}\label{thm:present}
	Let $\mathcal M$ be a countably presented MR-algebra. Then 
	$\mathcal M$ is a filter algebra. 
\end{thm}
\begin{proof}
	Let $A=\Set{a_{i} | i\in\omega}$ be a countable set that presents $\mathcal 
	M$. Define the sequence $\brk<b_{n} | n\in\omega>$ by 
	\begin{align*}
		b_{0} & =a_{0}  \\
		b_{n+1} & =b_{n}\caret a_{n+1}. 
	\end{align*}
	Then we have $b_{n+1}\le b_{n}$ and $b_{n+1}\preccurlyeq a_{n+1}$ so 
	that 
	$$
		\mathcal M=\bigcup_{n=0}^{\infty}\mathcal M_{a_{n}}=\bigcup_{n=0}^{\infty}\mathcal M_{b_{n}} 
	$$
	and $\Set{x | \exists n \ x\geq b_{n}}$ is therefore a generating 
	filter for $\mathcal M$. 
\end{proof}

\section{Some Facts about Filter Algebras}
Let $\rsf I(\rsf F)$ be a filter algebra. We recall some earlier 
results about automorphisms of $\rsf I(\rsf F)$ from \cite{BO:fil}. 

\begin{defn}[\cite{BO:fil}*{Definitions 22, 34}]\label{def:gensEtc}
	Let $\mathcal L=\rsf I(\rsf F)$ be any filter algebra. Let
	$\rsf G\subseteq\mathcal L$ be a  filter. 
	\begin{enumerate}[(a)]
		\item $\leftGen \rsf G\rightGen$ is the subalgebra generated by 
		$\rsf G$;
	
		\item $\widehat{\rsf G}=\Set{\Delta(x, y) | x, y\in\rsf G\text{ and 
		}y\le x}$; 
		
		\item $\rsf G$ is a \emph{generating filter} or \emph{g-filter} iff
		$\leftGen\rsf G\rightGen=\mathcal L$. 
	\end{enumerate}
\end{defn}

\begin{lem}[\cite{BO:fil}*{Theorem 29}]\label{lem:gen}
	Let $\rsf G$ be a filter in a cubic implication algebra $\mathcal L$. Then
	$$
	\leftGen\rsf G\rightGen=\widehat{\rsf G}. 
	$$
\end{lem}

\begin{defn}[{\cite{BO:fil}*{Definition 35}}]\label{def:AlphaBeta}
	Let $\mathcal M$ be an MR-algebra, and $\rsf F$ be a generating 
	filter in $\mathcal M$. Then for all $x\in\mathcal M$ we let
	$\alpha_{\rsf F}(x)$ and $\beta_{\rsf F}(x)$ be the unique 
	elements in $\rsf F$ such that
	$x=\Delta(\alpha_{\rsf F}(x), \beta_{\rsf F}(x))$. 
\end{defn}

\begin{thm}[{\cite{BO:fil}*{Corollary 43}}]\label{cor:filterAuts}
	Let $\mathcal M$ be an MR-algebra, and $\rsf F$, $\rsf G$ be two 
	generating filters in $\mathcal M$. Then
	the function $x\mapsto\beta_{\rsf G}(x)$ from $\rsf F$ to $\rsf G$
	is a one-one onto implication homomorphism. 
\end{thm}

\begin{defn}[{\cite{BO:fil}*{Definition 47}}]\label{def:filterAut}
	Let $\mathcal M$ be an MR-algebra, and $\rsf F$, $\rsf G$ be two 
	generating filters in $\mathcal M$. 
	
	A function $f\colon{\mathcal M}\to{\mathcal M}$ is a \emph{filter 
	automorphism based on $\brk<\rsf F, \rsf G>$} iff
	\begin{enumerate}[(a)]
		\item $f$ is an automorphism of $\mathcal M$;
	
		\item $f[\rsf F]=\rsf G$;
	
		\item for all $x\in\rsf F$\ $x\sim f(x)$. 
	\end{enumerate}
	Filter automorphisms are also called \emph{inner automorphisms} and 
	the set of all inner automorphisms is denoted by $\inn(\mathcal M)$. 
\end{defn}

\begin{lem}[\cite{BO:fil}*{Lemma 48, Definition 49}]\label{lem:uniqueness}
	Let $\mathcal M$ be an MR-algebra and $\rsf F$, $\rsf G$ be two 
	generating filters in $\mathcal M$. Then there is a unique filter automorphism 
	$f$ such that $f[\rsf F]=\rsf G$. This automorphism is denoted by 
	$\varphi_{\brk<\rsf F, \rsf G>}$. 
\end{lem}

\begin{lem}[{\cite{BO:fil}*{Lemma 51}}]\label{lem:anyFilter}
	Let $\mathcal M$ be a filter algebra. Let $\rsf F$ be a
	generating
	filter in $\mathcal M$. 
	Let $f$ be a cubic automorphism such that $x\sim f(x)$ for all 
	$x\in\mathcal M$. Then
	\begin{enumerate}[(a)]
		\item $f[\rsf F]$ is a generating filter; and
	
		\item $f=\varphi_{\brk<\rsf F, f[\rsf F]>}$. 
	\end{enumerate}
\end{lem}

From this lemma it is easy to show that the set of filter 
automorphisms is a group, but we also want to know that it has 
$2$-torsion. 
\begin{lem}[{\cite{BO:fil}*{Lemma 58, Corollary 59}}]\label{lem:inverses}
	Let $\mathcal M$ be an MR-algebra and $\rsf F, \rsf G$ be two 
	generating filters in $\mathcal M$. Then
	$$
	\varphi_{\brk<\rsf F, \rsf G>}=\varphi_{\brk<\rsf G, \rsf F>} 
	$$
	and hence 
	$$
	\varphi_{\brk<\rsf F, \rsf G>}^{-1}=\varphi_{\brk<\rsf F, \rsf G>}. 
	$$
\end{lem}

% The last results we need are about fixed points. 
% \begin{lem}[{\cite[Lemma 3.15]{BO:fil}}]\label{lem:above}
% 	Let $\mathcal M$ be an MR-algebra and $\rsf F, \rsf G$ be two 
% 	generating filters in $\mathcal M$. Then
% 	$$
% 	\varphi_{\brk<\rsf F, \rsf G>}\restrict(\rsf F\cap\rsf G)=\text{id}. 
% 	$$
% \end{lem}

Extending the ideas of \cite{BO:fil} we define operations and 
properties of filters and g-filters, and can prove certain 
consequences.  These results will appear in more detail in 
a later paper. 

\begin{defn}\label{def:impl}
	Let $\rsf G\subseteq\rsf F$ be two $\mathcal L$-filters. Then
	\begin{enumerate}[(a)]
		\item $\rsf G\supset\rsf F=\bigcap\Set{\rsf H | \rsf H\join\rsf G=\rsf F}$; 
	
		\item $\rsf G\Rightarrow\rsf F=\bigvee\Set{\rsf H |\rsf H\subseteq\rsf 
		F\text{ and }\rsf H\cap\rsf G={\Set{\one}}}$; 
	
		\item $\rsf G\to\rsf F=\Set{h\in\rsf F |\forall g\in\rsf G\ h\join 
		g=\one}$. 
	\end{enumerate}
\end{defn}

\begin{lem}\label{lem:twoThreeSame}
	\begin{enumerate}[(a)]
		\item[]
		\item $\rsf G\to\rsf F=\rsf G\Rightarrow\rsf F$. 
	
		\item $\rsf G\supset\rsf F=\rsf G\to\rsf F$. 
	\end{enumerate}	
\end{lem}

Particular amongst all filters are Boolean filters. 
\begin{defn}\label{def:Boolean}
	Let $\rsf F$ be a g-filter. Then
	\begin{enumerate}[(a)]
		\item $\rsf G$ is \emph{weakly $\rsf F$-Boolean} iff $\rsf G\subseteq\rsf 
		F$ and
		$(\rsf G\to\rsf F)\to\rsf F=\rsf G$. 
	
		\item $\rsf G$ is \emph{weakly Boolean} iff 
		there is some g-filter containing $\rsf G$ and 
		$\rsf G$ is $\rsf H$-Boolean for all such g-filters $\rsf H$.
		
		\item $\rsf G$ is \emph{$\rsf F$-Boolean} iff $\rsf G\subseteq\rsf 
		F$ and $\rsf G\join(\rsf G\to\rsf F)=\rsf F$. 
		
		\item $\rsf G$ is \emph{Boolean} iff
		there is some g-filter containing $\rsf G$ and 
		$\rsf G$ is $\rsf H$-Boolean for all such g-filters $\rsf H$.
	\end{enumerate}
\end{defn}

\begin{thm}\label{thm:Boolean}
	Let $\rsf G$ be $\rsf F$-Boolean for some g-filter $\rsf F$. 
	Then $\rsf G$ is Boolean. 
\end{thm}

Now we can define a $\Delta$ operation on filters and this allows 
recovery of a g-filter some certain fragments. 
\begin{defn}\label{def:Delta}
	Let $\rsf G\subseteq\rsf F$. Then 
	$$
	\Delta(\rsf G, \rsf F)=\Delta(\one, \rsf G\to\rsf F)\join\rsf G. 
	$$
\end{defn}

\begin{thm}\label{thm:lots}
	If $\rsf H$, $\rsf F$ are g-filters then 
	$\rsf G=\rsf F\cap\rsf H$ is $\rsf F$-Boolean and $\rsf H
	=\Delta(\rsf G, \rsf F)$.
	
	Conversely, if $\rsf G$ is $\rsf F$-Boolean then 
	$\rsf H=\Delta(\rsf G, \rsf F)$ is a g-filter and 
	$\rsf G=\rsf F\cap\rsf H$. 
\end{thm}

\section{Some Category Theory}
The operation $\rsf I$ is a functor where we define $\rsf I(f)\colon\rsf 
I(\mathcal I_{1})\to\rsf I(\mathcal I_{2})$ by
$$
		\rsf I(f)(\brk<a, b>)=\brk<f(a), f(b)>	
$$
whenever $f\colon\mathcal I_{1}\to\mathcal I_{2}$ is an implication 
morphism. 

The relation $\sim$ defined above
gives rise to a functor $\rsf C$ on cubic implication algebras. 
In order to show that it is well-defined this we need the following lemma. 
\begin{lem}\label{lem:simHom}
	Let $\phi\colon\mathcal L_{1}\to\mathcal L_{2}$ be a cubic 
	homomorphism. Let $a, b\in\mathcal L_{1}$. Then
	$$
		a\sim b\Rightarrow \phi(a)\sim\phi(b). 
	$$
\end{lem}
\begin{proof}
\begin{align*}
	a\sim b & \iff \Delta(a\join b, a)=b  \\
		 & \hphantom{\Leftarrow}\Rightarrow \phi(\Delta(a\join b, a))=\phi(b)  \\
		 & \iff \Delta(\phi(a)\join\phi(b), \phi(a))=\phi(b)  \\
		 & \iff \phi(a)\sim\phi(b). 
	\end{align*}
\end{proof}

$\rsf C$ is defined by
\begin{align*}
	\rsf C(\mathcal L) & =\mathcal L/\sim  \\
	\rsf C(\phi)([x]) & =[\phi(x)]. 
\end{align*}

There are several natural transformations here. The basic ones are 
$e\colon\text{ID}\to\rsf I$ and $\eta\colon\text{ID}\to\rsf C$, defined by 
\begin{align*}
	e_{\mathcal I}(x)&=\brk<\one, x>\\
	\eta_{\mathcal L}(x)&=[x]. 
\end{align*}
The following diagram commutes: 
\begin{diagram}
	\mathcal I_{1} & \rTo^{\phi} & \mathcal I_{2}  \\
	\dTo^{e_{\mathcal I_{1}}} &  &  \dTo^{e_{\mathcal I_{2}}} \\
	\rsf I(\mathcal I_{1}) & \rTo_{\rsf I(\phi)} & \rsf I(\mathcal I_{2}) 
\end{diagram}
as for $x\in\mathcal I_{1}$, we have 
\begin{align*}
	e_{\mathcal I_{2}}(\phi(x)) & =\brk<\one, \phi(x)>  \\
	 & =\brk<\phi(\one), \phi(x)>  \\
	 & =\rsf I(\phi)(\brk<\one, x>)\\
	 &=\rsf I(\phi)e_{\mathcal I_{1}}(x). 
\end{align*}

The following diagram commutes: 
\begin{diagram}
	\mathcal L_{1} & \rTo^{\phi} & \mathcal L_{2}  \\
	\dTo^{\eta_{\mathcal L_{1}}} &  &  \dTo^{\eta_{\mathcal L_{2}}} \\
	\rsf C(\mathcal L_{1}) & \rTo_{\rsf C(\phi)} & \rsf C(\mathcal L_{2}) 
\end{diagram}
as for $x\in\mathcal L_{1}$, we have
\begin{align*}
	\eta_{\mathcal L_{2}}(\phi(x)) & =[\phi(x)]  \\
	 & =\rsf C(\phi)([x])\\
	 &=\rsf C(\phi)\eta_{\mathcal L_{1}}(x). 
\end{align*}

Then we get the composite transformation 
$\iota\colon\text{ID}\to\rsf C\rsf I$ defined by 
$$
	\iota_{\mathcal I}=\eta_{\rsf I(\mathcal I)}\circ e_{\mathcal I}. 
$$
By standard theory this is a natural transformation. It is easy to 
see that $e_{\mathcal I}$ is an embedding, and that $\eta_{\mathcal L}$ 
is onto. 
\begin{thm}\label{thm:isoIota}
	$\iota_{\mathcal I}$ is an isomorphism. 
\end{thm}
\begin{proof}
	Let $x, y\in\mathcal I$ and suppose that $\iota(x)=\iota(y)$. 
	Then 
	\begin{align*}
		\iota(x) & =\eta_{\rsf I(\mathcal I)}(e_{\mathcal I}(x))  \\
		 & =[\brk<\one, x>]  \\
		 & =[\brk<\one, y>]. 
	\end{align*}
	Thus $\brk<\one, x>\sim\brk<\one, y>$. Now 
	\begin{align*}
		\Delta(\brk<\one, x>\join\brk<\one, y>, \brk<\one, y>) & =
		\Delta(\brk<\one, x\join y>, \brk<\one, y>)\\
		 & =\brk<(x\join y)\to y, x\join y>. 
	\end{align*}
	This equals $\brk<\one, x>$ iff $x=x\join y$ (so that $y\le x$) and 
	$(x\join y)\to y=\one$ so that $y=x\join y$ and $x\le y$. Thus $x=y$. 
	Hence $\iota$ is one-one. 
	
	It is also onto, as if $z\in\rsf C\rsf I(\mathcal I)$ then we have 
	$z=[w]$ for some $w\in\rsf I(\mathcal I)$. But we know that 
	$w=\brk<x, y>\sim\brk<\one, x\meet y>$ -- since 
	$\Delta(\brk<\one, y>, \brk<\one, x\meet y>)=\brk<x, y>$ -- and so
	$z=[\brk<\one, x\meet y>]=\eta_{\rsf I(\mathcal I)}(e_{\mathcal 
	I}(x\meet y))$. 
\end{proof}

We note that there is also a natural transformation 
$\kappa\colon\text{ID}\to\rsf I\rsf C$ defined by 
$$
	\kappa_{\mathcal L}=e_{\rsf C(\mathcal L)}\circ \eta_{\mathcal L}. 
$$
In general this is not an isomorphism as there are MR-algebras $\mathcal 
M$ that are not filter algebras, but $\rsf I(\rsf C(\mathcal M))$ is 
always a filter algebra. 

We also note that $\iota_{\rsf C(\mathcal L)}=\rsf C(\kappa_{\mathcal L})$ 
for all cubic implication algebras $\mathcal L$. 
The pair $\rsf I$ and $\rsf C$ do not form an adjoint pair. 

\section{Automorphisms} 
In \cite{BO:fil} we showed that the automorphism group of a 
filter algebra $\rsf I(\rsf F)$ decomposes into a group of \emph{inner} 
automorphisms and the group of implication automorphisms of $\rsf F$. 
In the special case of an interval algebra $\rsf I(B)$ the former 
group is isomorphic to  $\brk<B, 0, +>$. For arbirary $\rsf I(\rsf F)$
the group of inner automorphisms is 
isomorphic to a $2$-torsion group of subfilters of $\rsf F$. 

Now we want to look at the most general case of MR-algebras and 
determine some of the structure of the automorphism group. 

Let $\mathcal M$ be any MR-algebra. The functor $\rsf C$ induces a 
group homomorphism $\rsf C\colon\aut(\mathcal M)\to\aut(\rsf C(\mathcal M))$. 
We want to look at the kernel of $\rsf C$. 

The method is somewhat indirect and first we consider what $\rsf C$ 
does in the case that $\mathcal M$ is a filter algebra. 

\subsection{$\rsf C$ on Filter algebras}
The kernel of $\rsf C$ on a filter algebra is relatively easy to 
compute. 

\begin{thm}\label{thm:kerFilter}
	Let $\mathcal M$ be a filter algebra. Then 
	$$
		\ker(\rsf C)=\inn(\mathcal M). 
	$$
\end{thm}
\begin{proof}
	Let $\phi$ be any inner automorphism. Then we know that $x\sim\phi(x)$ 
	for all $x$ so that $[x]=[\phi(x)]=\rsf C(\phi)([x])$ for all $x$. 
	Thus $\phi\in\ker(\rsf C)$. 
	
	Conversely if $\phi\in\ker(\rsf C)$ then $x\sim\phi(x)$ for all $x$. 
	Then by \lemref{lem:anyFilter} we know that $\phi$ is an 
	inner automorphism. 
\end{proof}

We recall that if $\mathcal M$ is a filter algebra and $\rsf F$ is a 
g-filter then $\phi_{\rsf F}\colon\mathcal M\to\rsf I(\rsf F)$ defined 
by 
$$
	\phi_{\rsf F}(x)=\brk<\Delta(\one, x)\join\beta_{\rsf F}(x), x\join\beta_{\rsf F}(x)>
$$
is an isomorphism -- the $\rsf F$-presentation of $\mathcal M$. 

$\rsf C$ gives an isomorphism from $\rsf C(\phi_{\rsf F})\colon\rsf C(\mathcal 
M)\to \rsf C\rsf I(\rsf F)$. 
Note that $\rsf C(\phi_{\rsf F})^{-1}=\rsf C(\phi_{\rsf F}^{-1})$.

Putting this together with $\iota_{\rsf F}$ we have an isomorphism $\iota_{\rsf F}^{-1}\rsf 
C(\phi_{\rsf F})\colon
\rsf C(\mathcal M)\to \rsf F$. This induces an isomorphism of 
automorphism groups 
\begin{align*}
	\Xi\colon\aut(\rsf C(\mathcal M)) & \to\aut(\rsf F)\text{ given by }  \\
	\Xi(\alpha) & = \iota_{\rsf F}^{-1}\rsf C(\phi_{\rsf F})\alpha\rsf C(\phi_{\rsf F}^{-1})\iota_{\rsf F}. 
\end{align*}

\begin{lem}\label{lem:phiE}
	Let $x\in\rsf F$. Then $$
	\phi_{\rsf F}^{-1}e_{\rsf F}(x)=x. 
	$$
\end{lem}
\begin{proof}
	\begin{align*}
		\phi_{\rsf F}(x) & =\brk<\Delta(\one, x)\join\beta_{\rsf F}(x), x\join\beta_{\rsf F}(x)>  \\
		 & =\brk<\Delta(\one, x)\join x, x\join x>  \\
		 & =\brk<\one, x>\\
		 &= e_{\rsf F}(x). 
	\end{align*}
\end{proof}

From this information we are able to identify the image of $\rsf C$. 

\begin{thm}\label{thm:factoring}
	Let $\phi\in\aut(\mathcal M)$. Let $\phi=\varphi_{\brk<\rsf F, \rsf 
	G>}\circ\widehat\chi$ where $\chi\in\aut(\rsf F)$. Then
	$$
		\Xi\rsf C(\phi)=\chi. 
	$$
\end{thm}
\begin{proof}
	Let $x\in\rsf F$. 
	Then 
	\begin{align*}
		\rsf C(\phi_{\rsf F}^{-1})\iota_{\rsf F}(x) & =\rsf C(\phi_{\rsf 
		F}^{-1})\eta_{\rsf I(\rsf F)}e_{\rsf F}(x)  \\
		 & =\rsf C(\phi_{\rsf F}^{-1})([e_{\rsf F}(x)])  \\
		 & =[\phi_{\rsf F}^{-1}e_{\rsf F}(x)] \\
		 & =[x] &&\text{by the lemma. }  \\
		\rsf C(\phi)([x]) & =\rsf C(\varphi_{\brk<\rsf F, \rsf G>})\rsf C(\widehat 
		\chi)([x])  \\
		 & =\rsf C(\widehat\chi)([x]) &&\text{as }\varphi_{\brk<\rsf F, \rsf 
		 G>}\in\ker(\rsf C) \\
		 & =[\widehat\chi(x)]\\
		 &=[\chi(x)] && \text{ as }x\in\rsf F\\
		\iota_{\rsf F}^{-1}\rsf C(\phi_{\rsf F})([\chi(x)]) &= 
		\iota_{\rsf F}^{-1}[\phi_{\rsf F}(\chi(x))]\\
		&=\iota_{\rsf F}^{-1}[e_{\rsf F}(\chi(x))]\\
		&=\iota_{\rsf F}^{-1}\eta_{\rsf I(\rsf F)}e_{\rsf F}\chi(x)\\
		&=\iota_{\rsf F}^{-1}\iota_{\rsf F}\chi(x)\\
		&=\chi(x)
	\end{align*}
\end{proof}

% Thus we have the diagram 
% \begin{diagram}
% 	\rsf I(\rsf F) & \rTo^{\beta_{e_{\rsf F}[\rsf F]}} & e_{\rsf F}[\rsf F]  \\
% 	\dTo^{\eta_{\sim}} &  & \dTo_{e_{\rsf F}^{-1}}  \\
% 	\rsf I(\rsf F)/\sim & \rTo_{\iota_{\rsf F}^{-1}} & \rsf F
% \end{diagram}
% which commutes. 
% 
% Now let $\phi$ be any automorphism of $\mathcal M$. Then the diagram
% \begin{diagram}
% 	\mathcal M & \rTo^{\phi} & \mathcal M & \rTo^{\phi_{\rsf F}} & \rsf I(\rsf F) & \rTo^{\beta_{e_{\rsf F}[\rsf F]}} & e_{\rsf F}[\rsf F]  \\
% 	\dTo^{\eta_{\sim}} &  & \dTo^{\eta_{\sim}} &  & \dTo^{\eta_{\sim}} &  & \dTo_{e_{\rsf F}^{-1}} \\
% 	\mathcal M/\sim & \rTo_{\rsf C(\phi)} & \mathcal M/\sim & \rTo_{[\phi_{\rsf F}]} & \rsf I(\rsf F)/\sim & \rTo_{\iota_{\rsf F}^{-1}} & \rsf F
% \end{diagram}
% also commutes. 

Inner automorphisms of filter algebras are determined by their 
action on a single g-filter $\rsf F$. 

The results cited in section 3, in particular \thmref{thm:lots}, shows that 
any g-filter $\rsf G$ is determined by $\rsf F\cap\rsf G$. This set can 
be found from the set of fixed points for $\varphi_{\brk<\rsf F, \rsf G>}$. 
\begin{lem}\label{lem:fixed}
	Let $\mathcal M=\rsf I(\rsf F)$ be a filter algebra and 
	$\phi=\varphi_{\brk<\rsf F, \rsf G>}$ be any inner automorphism. Then
	$$
		\phi(x)=x\text{ iff }x\in\leftGen\rsf F\cap\rsf G\rightGen. 
	$$
\end{lem}
\begin{proof}
	Let $x\in\leftGen\rsf F\cap\rsf G\rightGen$. Then 
	$\beta_{\rsf F}(x)=\beta_{\rsf G}(x)$ and so 
	$\alpha_{\rsf F}(x)=\alpha_{\rsf G}(x)$. Thus
	\begin{align*}
		\phi(x) & =\Delta(\beta_{\rsf G}\alpha_{\rsf F}(x), \beta_{\rsf G}(x))  \\
		 & =\Delta(\alpha_{\rsf F}(x), \beta_{\rsf F}(x))  \\
		 & =x. 
	\end{align*}
	Conversely,  if $\phi(x)=x$ then $x= \phi(x)= \Delta(\beta_{\rsf G}\alpha_{\rsf F}(x), \beta_{\rsf G}(x))$
	so that $\beta_{\rsf G}\alpha_{\rsf F}(x)=\alpha_{\rsf G}(x)$. Since 
	$x\le \alpha_{\rsf F}(x)\meet\alpha_{\rsf G}(x)$ this implies 
	$\alpha_{\rsf F}(x)=\alpha_{\rsf G}(x)$ is in $\rsf F\cap\rsf G$. 
	But now we have $\beta_{\rsf F}(x)= \Delta(\alpha_{\rsf F}(x), x)= 
	\Delta(\alpha_{\rsf G}(x), x)= \beta_{\rsf G}(x)$ and so $x\in\leftGen\rsf F\cap\rsf G\rightGen$. 
\end{proof}

\begin{lem}\label{lem:DeltaFixed}
	Let $\mathcal M=\rsf I(\rsf F)$ be a filter algebra and 
	$\phi=\varphi_{\brk<\rsf F, \rsf G>}$ be any inner automorphism. Then
	$$
		\phi(x)=\Delta(\one, x)\text{ iff }x\in\leftGen(\rsf F\cap\rsf G)\to\rsf F\rightGen. 
	$$
\end{lem}
\begin{proof}
	First we recall that $\Delta(\one, \rsf G)\cap\rsf F=(\rsf F\cap\rsf 
	G)\to\rsf F$ -- see \cite[lemma 5.26]{BO:fil}. 
	
	Let $x\in(\rsf F\cap\rsf G)\to\rsf F=\Delta(\one, \rsf G)\cap\rsf F$. Then 
	$\phi(x)=\beta_{\rsf G}(x)$. As $\Delta(\one, x)\in\rsf G$ we have 
	$\beta_{\rsf G}(x)=\Delta(\one, x)$. 
	
	In general we have $x\in\leftGen(\rsf F\cap\rsf G)\to\rsf F\rightGen$ implies 
	$x=\Delta(\alpha_{\rsf F}(x), \beta_{\rsf F}(x))$ where both 
	$\alpha_{\rsf F}(x)$ and $\beta_{\rsf F}(x)$ are in 
	$(\rsf F\cap\rsf G)\to\rsf F$. Then we have 
	\begin{align*}
		\phi(x) & = \Delta(\beta_{\rsf G}\alpha_{\rsf F}(x), \beta_{\rsf G}\beta_{\rsf F}(x)) \\
		 & =\Delta(\Delta(\one, \alpha_{\rsf F}(x)), \Delta(\one, \beta_{\rsf F}(x)))  \\
		 & =\Delta(\one, \Delta(\alpha_{\rsf F}(x), \beta_{\rsf F}(x)))\\
		 &= \Delta(\one, x). 
	\end{align*}
	
	Conversely, if $\phi(x)=\Delta(\one, x)$ then we have 
	$\Delta(\one, x)= \phi(x)= \Delta(\beta_{\rsf G}\alpha_{\rsf F}(x), \beta_{\rsf G}(x))$
	so that $\beta_{\rsf G}\alpha_{\rsf F}(x)=\alpha_{\rsf G}(\Delta(\one, x))$. Since 
	$\beta_{\rsf G}\alpha_{\rsf F}(x)\sim\Delta(\one, \alpha_{\rsf 
	F}(x))$ and $\Delta(\one, x)\le \alpha_{\rsf G}(\Delta(\one, x))\meet \Delta(\one, 
	\alpha_{\rsf F}(x))$ this gives $\alpha_{\rsf G}(\Delta(\one, x))=\Delta(\one, 
	\alpha_{\rsf F}(x))$. Now we have 
	$x= \Delta(\alpha_{\rsf F}(x), \Delta(\one, \beta_{\rsf G}(x)))$ so 
	that 
	$\beta_{\rsf F}(x)= \Delta(\alpha_{\rsf F}(x), x)= \Delta(\one, 
	\beta_{\rsf G}(x))$. Thus
	$x\in\leftGen\Delta(\one, \rsf G)\cap\rsf F\rightGen$. 
\end{proof}

\subsection{Localizing}
The main technique we will use to derive information about 
automorphism is localization. Previously we have localized at a point 
to obtain interval algebras of the form $\mathcal L_{x}$ that are 
upwards-closed and contain $x$. Now we need more closure, so we will 
use a localization that produces filter algebras. 

The essential use of localization is contained in the next two results. 
\begin{thm}\label{thm:incl}
	Let $\mathcal L$ be a cubic implication algebra and $\mathcal M$ be any 
	upwards-closed subalgebra. 
	Then $\rsf C(\text{incl}_{\mathcal M})=\text{incl}_{\rsf C(\mathcal 
	M)}$ and the following diagram commutes:
	\begin{diagram}
		\mathcal M & \rTo^{\incl_{\mathcal M}} & \mathcal L  \\
		\dTo^{\eta_{\mathcal M}} &  & \dTo^{\eta_{\mathcal L}}  \\
		\rsf C(\mathcal M) & \rTo_{\incl_{\rsf C(\mathcal M)}} & \rsf C(\mathcal L)
	\end{diagram}
\end{thm}
\begin{proof}
	The commutativity of the diagram follows from the first result as 
	$\eta$ is a natural transformation from $\text{ID}$ to $\rsf C$. 
	
	We show that if $x\in\mathcal M$ then $[x]_{\mathcal M}=[x]_{\mathcal L}$. 
	
	Let $y\in[x]_{\mathcal M}$. Then $x\sim y$ in $\mathcal M$. But this 
	happens iff $\Delta(x\join y, y)=x$ which is true in $\mathcal M$ 
	iff it is true in $\mathcal L$. Thus $y\in[x]_{\mathcal L}$. 
	
	Let $y\in[x]_{\mathcal L}$. Then $x\join y\in\mathcal M$ as $\mathcal 
	M$ is upwards-closed. But then $\Delta(x\join y, x)=y$ is also in $\mathcal 
	M$ and so $y\in[x]_{\mathcal M}$. 
\end{proof}

\begin{cor}\label{cor:restrict}
	Let $f\colon\mathcal L_{1}\to\mathcal L_{2}$ be a cubic homomorphism. 
	Let $\mathcal M$ be any upwards-closed subalgebra of $\mathcal L_{1}$. 
	Then 
	$$
		\rsf C(f\restrict\mathcal M)=\rsf C(f)\restrict\rsf C(\mathcal M). 
	$$
\end{cor}
\begin{proof}
	$f\restrict\mathcal M=f\circ\text{incl}_{\mathcal M}$ so that 
	\begin{align*}
		\rsf C(f\restrict\mathcal M) & =\rsf C(f\circ\text{incl}_{\mathcal M})  \\
		 & =\rsf C(f)\circ\rsf C(\text{incl}_{\mathcal M})  \\
		 & =\rsf C(f)\circ\text{incl}_{\rsf C(\mathcal M)} \\
		 & =\rsf C(f)\restrict\rsf C(\mathcal M). 
	\end{align*}
\end{proof}

Now here is the localization process we need. 
\begin{thm}\label{thm:localization}
	Let $\mathcal M$ be any MR-algebra. Let $X$ be a countable subset of 
	$\mathcal M$ and $\mathcal G$ be any countable subgroup of $\aut(\mathcal 
	M)$. Then there is a subalgebra $\mathcal L$ of $\mathcal M$ such that
	\begin{enumerate}[(a)]
		\item $\mathcal L$ is an upwards-closed MR-subalgebra that is 
		countably presented; 
	
		\item $X\subseteq\mathcal L$; 
	
		\item if $\phi\in\mathcal G$ then $\phi\restrict\mathcal L$ is in 
		$\aut(\mathcal L)$. 
	\end{enumerate}
\end{thm}
\begin{proof}
	Define $Z$ inductively by:
	\begin{align*}
		Z_{0} & =X  \\
		Z_{2n+1} & =\text{the caret-closure of }Z_{2n}  \\
		Z_{2n+2} & = \Set{\phi^{i}(y) | y\in Z_{2n+1},\ \phi\in\mathcal G,\ i\in\omega} \\
		Z & =\bigcup_{m\in\omega}Z_{m}. \\
		\intertext{Now let }
		\mathcal L&=\bigcup_{z\in Z}\mathcal M_{z}.
	\end{align*}
	Then clearly $\mathcal L$ is upwards-closed. It is easy to see that 
	$X\subseteq Z_{m}\subseteq Z$ for all $m$ so that $X\subseteq \mathcal L$. Also $Z$ 
	is countable so that $\mathcal L$ is countably presented. 
	
	$\mathcal L$ is caret-closed -- as if $x, y\in\mathcal L$ then let 
	$a\preccurlyeq x$ and $b\preccurlyeq y$ for some $a, b\in Z_{m}$. Then 
	$m$ odd implies $a\caret b\in Z_{m}$, else $a\caret b\in Z_{m+1}$. 
	Thus $a\caret b\preccurlyeq x\caret y$ and so $x\caret y\in \mathcal L$. 
	
	Let $\phi\in\mathcal G$. Then $\mathcal L$ is $\phi$-closed -- as if 
	$a\in Z_{m}$ and $a\preccurlyeq z$ then either $\phi(a)\in Z_{m}$ (if 
	$m$ is even) or $\phi(a)\in Z_{m+1}$. Now we have 
	$\phi(a)\preccurlyeq\phi(z)$ so that $\phi(z)\in\mathcal L$. 
	
	Let $\phi\in\mathcal G$. Then $\phi\restrict\mathcal L$ is clearly a one-one 
	homomorphism from $\mathcal L$ to $\mathcal L$. As $\phi^{-1}\in\mathcal 
	G$ we know that $\phi\restrict\mathcal L$ is also onto. Thus 
	$\phi\restrict\mathcal L$ is in $\aut(\mathcal L)$. 
\end{proof}

\begin{cor}\label{cor:local}
	Let $\mathcal M$ be any MR-algebra. Let $X$ be a countable subset of 
	$\mathcal M$ and $\mathcal X$ be any countable subset of $\aut(\mathcal 
	M)$. Then there is a subalgebra $\mathcal L$ of $\mathcal M$ such that
	\begin{enumerate}[(a)]
		\item $\mathcal L$ is an upwards-closed MR-subalgebra that is 
		countably presented; 
	
		\item $X\subseteq\mathcal L$; 
	
		\item if $\phi\in\mathcal X$ then $\phi\restrict\mathcal L$ is in 
		$\aut(\mathcal L)$. 
	\end{enumerate}
\end{cor}
\begin{proof}
	Let $\mathcal G$ be the subgroup of $\aut(\mathcal M)$ generated by $\mathcal 
	X$ and apply the theorem. 
\end{proof}

\subsection{Inner Automorphisms}
We wish to analyze the kernel of $\rsf C$. Generalizing from filter algebras
we define the group of inner automorphisms. 
\begin{defn}\label{def:innerAut}
	The group $\ker(\rsf C)$ is called the group of \emph{inner 
	automorphisms} denoted by $\inn(\mathcal M)$. 
\end{defn}

Earlier (\lemref{lem:fixed} and above) we showed that inner 
automorphisms on filter algebras are determined by their set of fixed 
points. We will show that this is always true. 
\begin{defn}\label{def:}
	Let $\phi\in\inn(\mathcal M)$. Let
	$$
		\mathcal M_{\phi}=\Set{x | \phi(x)=x}. 
	$$
\end{defn}

\subsection{$2$-torsion}
Let $\mathcal M$ be any MR-algebra and 
let $\alpha\in\aut(\mathcal M)$ be an inner automorphism. We wish to 
show that $\alpha^{2}=\text{id}$. 

\begin{thm}\label{thm:TwoTorsion}
	Let $\alpha$ be any inner automorphism of $\mathcal M$. Then 
	$\alpha^{2}=\text{id}$. 
\end{thm}
\begin{proof}
	Let $x\in\mathcal M$. We show that $\alpha^{2}(x)=x$. 
	
	Let $\mathcal L$ be as provided by \corref{cor:local} for $x$ and $\alpha$. 
	Then since we have (by \corref{cor:restrict}) 
	$\rsf C(\alpha\restrict\mathcal L)=\rsf C(\alpha)\restrict\rsf C(\mathcal 
	L)$ and $\rsf C(\alpha)=\text{id}$ we know that $\alpha\restrict\mathcal 
	L$ is an inner automorphism of $\mathcal L$. So let $\rsf F$ and $\rsf 
	G$ be two $\mathcal L$-g-filters such $\alpha\restrict\mathcal 
	L=\varphi_{\brk<\rsf F, \rsf G>}$. Then we have 
	$$
		\alpha^{2}(x)=\varphi_{\brk<\rsf F, \rsf G>}^{2}(x)=x. 
	$$
\end{proof}

Thus $\inn(\mathcal M)$ is an abelian group normal in $\aut(\mathcal M)$. 

\section{The Group of Inner Automorphisms}
The technique of localization, as used above, establishes more 
about the group of inner automorphisms. Here we will use it to show 
how to recover $\phi$ from $\mathcal M_{\phi}$ and hence that 
$\inn(\mathcal M)\simeq\inn(\rsf I\rsf C(\mathcal M))$.

\begin{lem}\label{lem:upper}
	Let $\phi\in\inn(\mathcal M)$. Then
	$\mathcal M_{\phi}$ is an upwards-closed MR-subalgebra of $\mathcal M$. 
\end{lem}
\begin{proof}
	Let $x\in\mathcal M_{\phi}$ and let $\mathcal L$ be as given by 
	\corref{cor:local}  for $x$ and $\phi$. Then (as in \thmref{thm:TwoTorsion})
	we know that $\phi\restrict\mathcal L$ is an inner automorphism of $\mathcal 
	L$ --
	say it equals $\varphi_{\brk<\rsf F, \rsf G>}$. Then 
	\begin{align*}
		\mathcal M_{\phi}\cap\mathcal L & =\Set{x\in\mathcal L | \phi(x)=x}  \\
		 & =\leftGen\rsf F\cap\rsf G\rightGen && \text{ by \lemref{lem:fixed}}
	\end{align*}
	which is upwards-closed. Thus $[x, \one]\subseteq\mathcal M_{\phi}$. 
\end{proof}

The sets $\mathcal M_{\phi}$ are non-trivial. For example, if $x\in\mathcal 
M$ then $x\join\phi(x)$ is in $\mathcal M_{\phi}$ -- since 
$\phi^{2}=\text{id}$. Also $x\join\Delta(\one, \phi(x))$ is 
not in $\mathcal M_{\phi}$ unless $x=\one$.

\begin{thm}\label{thm:MPhiIsGood}
	Let $\phi_{1}$ and $\phi_{2}$ be in $\inn(\mathcal M)$. Then
	$$
		\mathcal M_{\phi_{1}}=\mathcal M_{\phi_{2}}\text{ implies 
		}\phi_{1}=\phi_{2}. 
	$$
\end{thm}
\begin{proof}
	Let $x\in\mathcal M$. Let $\mathcal L$ be as given by 
	\corref{cor:local} for $x$ and $\Set{\phi_{1}, \phi_{2}}$. Let $\rsf 
	F$, $\rsf G_{1}$ and $\rsf G_{2}$ be $\mathcal L$-g-filters such that
	$\phi_{i}\restrict\mathcal L=\varphi_{\brk<\rsf F, \rsf G_{i}>}$. 
	Then we have  
	$\mathcal M_{\phi_{i}}\cap\mathcal L=\leftGen\rsf F\cap\rsf 
	G_{i}\rightGen$ and as $\mathcal M_{\phi_{1}}=\mathcal M_{\phi_{2}}$ 
	we have 
	$\rsf F\cap\rsf G_{1}=\rsf F\cap\rsf G_{2}$. But then 
	$\rsf G_{1}= \Delta(\rsf F\cap\rsf G_{1}, \rsf F)= \Delta(\rsf F\cap\rsf 
	G_{2}, \rsf F)= \rsf G_{2}$ -- see\thmref{thm:lots}. 
	
	Hence $\phi_{1}\restrict\mathcal L=\phi_{2}\restrict\mathcal L$ and so
	$\phi_{1}(x)=\phi_{2}(x)$. 
	
	Since we can do this for any $x\in\mathcal M$ we have $\phi_{1}=\phi_{2}$. 
\end{proof}

This theorem gives only a suggestion  of how to recover $\phi$ 
from $\mathcal M_{\phi}$. We'll now show how to fully recover $\phi$, 
as a prelude to showing that $\inn(\mathcal M)\simeq\inn(\rsf I\rsf 
C(\mathcal M))$. 

First we note that upwards closed subalgebras are completely 
determined by their collapses. 
\begin{lem}\label{lem:collapseDewt}
	Let $\mathcal L_{1}$ and $\mathcal L_{2}$ be two upwards-closed 
	subalgebras of a cubic implication algebra $\mathcal M$. Then
	$$
		\mathcal L_{1}=\mathcal L_{2}\text{ iff }\rsf C(\mathcal L_{1})=\rsf 
		C(\mathcal L_{2}). 
	$$
\end{lem}
\begin{proof}
	This is because $\mathcal L_{i}=\bigcup\Set{[x] | [x]\in\rsf C(\mathcal 
	L_{i})}$. 
\end{proof}

The main fact we know about $\phi$ is that $x\sim\phi(x)$ for all 
$x$. This is another way of viewing
$\rsf C(\phi)=\text{id}$. It implies that 
\begin{align}
	 x&=(x\join\phi(x))\meet(x\join\Delta(\one, \phi(x))) 
	\label{eq:oneA}  \\
	\phi(x)&=(x\join\phi(x))\meet(\Delta(\one, x)\join\phi(x)) 
	\label{eq:two}  \\
	\Delta(\one, x)\join\phi(x)&=\Delta(\one, x\join\Delta(\one, \phi(x))) 
	\label{eq:three}
\end{align}
for any $x$. Since we already know that $x\join\phi(x)$ is in $\mathcal M_{\phi}$
we need to look at the other half of the first equation. 
The other two equations suggest the recovery of $\phi$ -- it's the 
identity on $\mathcal M_{\phi}$ and $\Delta(\one, \bullet)$ on the 
other side. So we need to identify the other side!

\begin{defn}\label{def:otherSide}
	Let $\phi\in\inn(\mathcal M)$. Let 
	$$
		D_{\phi}=\eta_{\mathcal M}^{-1}[\rsf C(\mathcal M_{\phi})\to\rsf C(\mathcal M)]. 
	$$
\end{defn}

\begin{lem}\label{lem:somethingIn}
	Let $x\in\mathcal M$ be arbitrary. Then $\Delta(\one, 
	x)\join\phi(x)\in D_{\phi}$. 
\end{lem}
\begin{proof}
	Let $z=\eta_{\mathcal M}(\Delta(\one, x)\join\phi(x))$. We want to 
	show that $z\in\rsf C(\mathcal M_{\phi})\to\rsf C(\mathcal M)$ -- 
	that is for any $[y]\in\rsf C(\mathcal M_{\phi})$ we have 
	$z\join[y]=\one$. 
	This is equivalent to showing that $\eta_{\mathcal M}(y*(\Delta(\one, 
	x)\join\phi(x)))=\one$ and as $\eta_{\mathcal 
	M}^{-1}[\one]=\Set{\one}$ we need to show that 
	$y*(\Delta(\one, x)\join\phi(x))=\one$. 
	
	\begin{align*}
		y*(\Delta(\one, x)\join\phi(x)) & =\Delta(\Delta(\one, 
		x)\join\phi(x)\join y, y)\join (\Delta(\one, x)\join\phi(x))  \\
		\Delta(\one, 
		x)\join\phi(x)\join y & \geq y  \\
		\intertext{so it must be in $\mathcal M_{\phi}$ and therefore}
		\Delta(\one, x)\join\phi(x)\join y & =\phi(\Delta(\one, 
		x)\join\phi(x)\join y)  \\
		 & =\Delta(\one, \phi(x))\join x\join y  \\
		\text{ and so } 
		\Delta(\one, x)\join\phi(x)\join y & =(\Delta(\one, x)\join\phi(x)\join 
		y)\join (\Delta(\one, \phi(x))\join x\join y)   \\
		 & =\one. \\
		 \intertext{Thus }
		 y*(\Delta(\one, x)\join\phi(x))&=\Delta(\one, y)\join
		 \Delta(\one, x)\join\phi(x). 
	\end{align*}
	As this is also in $\mathcal M_{\phi}$ it must also equal $\one$. 
\end{proof}

\begin{lem}\label{lem:DeltaOne}
	Let $z\in D_{\phi}$. Then $\phi(z)=\Delta(\one, z)$. 
\end{lem}
\begin{proof}
	Let $\mathcal L$ be as given by \corref{cor:local} for $z$ and $\phi$. 
	Let $\phi\restrict\mathcal L=\varphi_{\brk<\rsf F, \rsf G>}$. We 
	know that $\mathcal M_{\phi}\cap\mathcal L=\leftGen\rsf F\cap\rsf 
	G\rightGen$ so that -- in $\mathcal L$ -- we have 
	$D_{\phi\restrict\mathcal L}=\leftGen(\rsf F\cap\rsf G)\to\rsf 
	F\rightGen$ and that 
	$\phi\restrict\mathcal L$ is equal to $\Delta(\one, \bullet)$ on this 
	set -- by \lemref{lem:DeltaFixed}. 
	
	Since $z\in D_{\phi}$ we have $\eta_{\mathcal M}(z)=\eta_{\mathcal 
	L}(z)\in \rsf C(\mathcal M_{\phi}\cap\mathcal L)\to\rsf C(\mathcal L)$ 
	and so $z\in D_{\phi\restrict\mathcal L}$. 
	Thus $\phi(z)= \left(\phi\restrict\mathcal L\right)(z)= \Delta(\one, z)$. 
\end{proof}

\begin{cor}\label{cor:intersect}
	$\mathcal M_{\phi}\cap D_{\phi}=\Set{\one}$. 
\end{cor}
\begin{proof}
	Let $x\in \mathcal M_{\phi}\cap D_{\phi}$. Then we have 
	$x= \phi(x)= \Delta(\one, x)$ so that 
	$x= x\join\Delta(\one, x)= \one$. 
\end{proof}

\begin{cor}\label{cor:metsExist}
	If $x\in\mathcal M_{\phi}$ and $y\in D_{\phi}$ then $x\meet 
	\Delta(\one, y)$ exists. 
\end{cor}
\begin{proof}
	This follows from the MR-axiom as $x\join y\in \mathcal 
	M_{\phi}\cap D_{\phi}$ and so equals $\one$. 
\end{proof}

\begin{lem}\label{lem:repsMD}
	Let $z\in\mathcal M$ be arbitrary. Then there is a unique pair
	$\brk<z_{0}, z_{1}>\in\mathcal M_{\phi}\times D_{\phi}$ such that
	$$
		z=z_{0}\meet z_{1}. 
	$$
\end{lem}
\begin{proof}
	We know that $z\sim\phi(z)$ so that 
	$z=(z\join\phi(z))\meet(z\join\Delta(\one, \phi(z)))$. From above 
	we have $z\join\phi(z)\in \mathcal M_{\phi}$ and 
	$\Delta(\one, z)\join\phi(z)\in D_{\phi}$. As $D_{\phi}$ is 
	$\Delta(\one, \bullet)$-closed we have  
	$z\join\Delta(\one, \phi(z))\in D_{\phi}$. 
	
	Suppose that $z=z_{0}\meet z_{1}$ with $\brk<z_{0}, z_{1}>\in
	\mathcal M_{\phi}\times D_{\phi}$. Then we have 
	\begin{align*}
		z\join\phi(z) & =(z_{0}\meet z_{1})\join(\phi(z_{0})\meet\phi(z_{1}))  \\
		 & = (z_{0}\meet z_{1})\join(z_{0}\meet \Delta(\one, z_{1}))  \\
		 & =z_{0}\meet( z_{1}\join\Delta(\one, z_{1}))  \\
		 & =z_{0}. 
	\end{align*}
	Likewise we have $z\join\Delta(\one, \phi(z))=z_{1}$. 
\end{proof}

\begin{lem}\label{lem:gotIt}
	Let $z\in\mathcal M$ be arbitrary and
	$\brk<z_{0}, z_{1}>\in\mathcal M_{\phi}\times D_{\phi}$ be such that
	\begin{align*}
		z&=z_{0}\meet z_{1}. \\
	\intertext{\qquad Then }
	\phi(z)&=z_{0}\meet \Delta(\one, z_{1}). 
	\end{align*}
\end{lem}
\begin{proof}
	Since $\phi$ preserves meets that exist and from the definition of $\mathcal 
	M_{\phi}$ and \lemref{lem:DeltaOne}. 
\end{proof}

This completes our recovery from $\phi$ from $\mathcal M_{\phi}$. But 
we need more in order to prove that 
$$
	\inn(\mathcal M)\simeq\inn(\rsf I\rsf C(\mathcal M)). 
$$

\begin{lem}\label{lem:BoolCC}
	$\rsf C(\mathcal M_{\phi})$ is $\rsf C(\mathcal M)$-Boolean. 
\end{lem}
\begin{proof}
	Let $z\in\mathcal M$. Then 
	\begin{align*}
		z& =(z\join\phi(z))\meet(z\join\Delta(\one, \phi(z)))  \\
		 & =(z\join\phi(z))\caret (\Delta(\one, z)\join\phi(z)). 
	\end{align*}
	Thus we have $[z]=[z\join\phi(z)]\meet[\Delta(\one, z)\join\phi(z)]$. 
	
	We know that $z\join\phi(z)\in\mathcal M_{\phi}$ and so 
	$[z\join\phi(z)]\in\rsf C(\mathcal M_{\phi})$. Also 
	$\Delta(\one, z)\join\phi(z)\in D_{\phi}$ so that 
	$[\Delta(\one, z)\join\phi(z)]\in \rsf C(\mathcal M_{\phi})\to\rsf C(\mathcal M)$. 
\end{proof}

\begin{lem}\label{lem:localBoolean}
	Let $\rsf F$ be a g-filter in some filter algebra $\mathcal L$, and 
	$\rsf G$, $\rsf H$ be two subfilters of $\rsf F$. 
	If $\rsf G$ is $\rsf F$-Boolean then $\rsf G\cap\rsf H$ is $\rsf 
	H$-Boolean. 
\end{lem}
\begin{proof}
	Firstly we see that if $s\in(\rsf G\to\rsf F)\cap\rsf H$ then for 
	all $g\in\rsf G\cap\rsf H$ we have $s\join g=\one$. Thus
	$$
		(\rsf G\to\rsf F)\cap\rsf H\subseteq(\rsf G\cap\rsf H)\to\rsf H. 
	$$
	Now if $h\in\rsf H$ then there is some $g\in\rsf G$ and $k\in\rsf 
	G\to\rsf F$ such that $h=g\meet k$. As $h\le g$ we have $g\in\rsf 
	G\cap\rsf H$. And $h\le k$ implies $k\in (\rsf G\to\rsf F)\cap\rsf 
	H$ and so $k\in(\rsf G\cap\rsf H)\to\rsf H$. 
\end{proof}

\begin{lem}\label{lem:localPrincBool}
	Let $\rsf F$ be a g-filter in some filter algebra $\mathcal L$, and 
	$\rsf G$ be a $\rsf F$-Boolean subfilter of $\rsf F$. Let $f\in\rsf F$. 
	Then $\rsf G\cap[f, \one]$ is principal. 
\end{lem}
\begin{proof}
	Let $f=g\meet h$ for $g\in\rsf G$ and $h\in\rsf G\to\rsf F$. Then we 
	have $[g, \one]\subseteq\rsf G\cap[f, \one]$. 
	
	If $x\in\rsf G\cap[f, \one]$ then 
	\begin{align*}
		x & =x\join f  \\
		 & =x\join(g\meet h)  \\
		 & =(x\join g)\meet (x\join h)  \\
		 & =x\join g &&\text{ as }x\join h=\one. 
	\end{align*}
	Thus $x\geq g$, and so $[g, \one]=\rsf G\cap[f, \one]$.
\end{proof}

\begin{lem}\label{lem:intComp}
	Let $\mathcal L=\mathcal L_{a}$ be an interval algebra. Let $g\geq 
	a$ and $h=g\to a$. Then for any $z\in\mathcal L$ 
	$$
		z=(z\join\Delta(g\join z, g))\meet(z\join\Delta(h\join z, h)). 
	$$
\end{lem}
\begin{proof}
	We work in an interval algebra $\rsf I(B)$. Wolog $a=[0, 0]$ so that
	$g=[0, g]$ and $h=[0, \comp g]$. Let $z=[z_{0}, z_{1}]$. 
	Then we have 
	\begin{align*}
		g\join z & =[0, g\join z_{1}]  \\
		h\join z & =[0, \comp g\join z_{1}]  \\
		\Delta(g\join z, g) & =[0\join((g\join z_{1})\meet\comp g), 
		0\join((g\join z_{1})\meet\comp 0)]  \\
		 & =[\comp g\meet z_{1}, g\join z_{1}]  \\
		 \Delta(h\join z, h) & =[g\meet z_{1},\comp  g\join z_{1}]   \\
		z\join\Delta(g\join z, g) & =[\comp g\meet z_{0}, g\join z_{1}]  \\
		z\join\Delta(h\join z, h) & =[g\meet z_{0},\comp  g\join z_{1}]  \\
		(z\join\Delta(g\join z, g))\meet(z\join\Delta(h\join z, h)) & =
		[(\comp g\meet z_{0})\join(g\meet z_{0}), (g\join z_{1})\meet(\comp  g\join z_{1})]\\
		 & =[z_{0}, z_{1}]. 
	\end{align*}
\end{proof}

\begin{thm}\label{thm:recoveryII}
	Let $\rsf G$ be a $\rsf C(\mathcal M)$-Boolean filter. Let
	\begin{align*}
		S_{1}&=\eta_{\mathcal M}^{-1}[\rsf G], \qquad \\
		S_{2}&=\eta_{\mathcal M}^{-1}[\rsf G\to\rsf C(\mathcal M)]. 
	\end{align*}
	Then 
	\begin{enumerate}[1.]
		\item $S_{1}\cap S_{2}=\Set{\one}$. 
		
		\item For all $x\in\mathcal M$ there are unique $x_{1}\in S_{1}$ 
		and $x_{2}\in S_{2}$ with $x=x_{1}\meet x_{2}$. 
	
		\item If we define $\phi_{\rsf G}\colon\mathcal M\to \mathcal M$ by
		\begin{equation}
			\phi_{\rsf G}(x)=x_{1}\meet \Delta(\one, x_{2})
			\label{eq:defPhiG}
		\end{equation}
		then 
		\begin{enumerate}[(a)]
			\item $\phi_{\rsf G}$ is a well-defined cubic automorphism of $\mathcal 
			M$; 
		
			\item $\phi_{\rsf G}$ is an inner automorphism of $\mathcal M$; 
		
			\item $\mathcal M_{\phi_{\rsf G}}=S_{1}$. 
		\end{enumerate}
	\end{enumerate}
\end{thm}
\begin{proof}
	\begin{enumerate}[B1.]
		\item Let $x\in S_{1}\cap S_{2}$. Then we have 
		$[x]=[x]\join[x]=\one$ so that $x\in\eta_{\mathcal 
		M}^{-1}[\one]=\Set{\one}$. 
		
		Note that this implies $x\meet\Delta(\one, y)$ exists for all $x\in 
		S_{1}$ and $y\in S_{2}$ from the MR-axiom and $x\join y=\one$. 
		
		\item Let $x\in\mathcal M$. Then $[x]=\mathbf{x}_{1}\meet \mathbf{x}_{2}$ for some 
		$\mathbf{x}_{1}\in\rsf G$ and $\mathbf{x}_{2}\in\rsf G\to\rsf C(\mathcal M)$ -- 
		as $\rsf G$ is a $\rsf C(\mathcal M)$-Boolean filter. 
		
		Let $x'_{1}\in S_{1}$ and $x'_{2}\in S_{2}$ be such that 
		$\mathbf{x}_{i}=[x'_{i}]$. Then we have $x\sim x'_{1}\caret x'_{2}$. 
		Hence 
		\begin{align*}
			x& =\Delta(x\join (x'_{1}\caret x'_{2}), x'_{1}\caret x'_{2})   \\
			 & =\Delta(x\join (x'_{1}\caret x'_{2}), x'_{1}\meet 
			 \Delta(x'_{1}\join x'_{2}, x'_{2}))  \\
			 & = \Delta(x\join (x'_{1}\caret x'_{2}), x'_{1})\meet 
			 \Delta(x\join (x'_{1}\caret x'_{2}), \Delta(x'_{1}\join x'_{2}, x'_{2})). 
		\end{align*}
		As $\Delta(x\join (x'_{1}\caret x'_{2}), x'_{1})\sim x'_{1}$ we 
		have $\Delta(x\join (x'_{1}\caret x'_{2}), x'_{1})\in S_{1}$. 
		Likewise $\Delta(x\join (x'_{1}\caret x'_{2}), \Delta(x'_{1}\join 
		x'_{2}, x'_{2}))$ is in $S_{2}$. 
		
		Suppose that $x= x_{1}\meet x_{2}= y_{1}\meet y_{2}$ with $x_{1}, 
		y_{1}$ in $S_{1}$ and $x_{2}, y_{2}$ in $S_{2}$. 
		Then we have 
		\begin{align*}
			x_{1} & =x_{1}\join x  \\
			 & =x_{1}\join(y_{1}\meet y_{2})  \\
			 &  =(x_{1}\join y_{1})\meet (x_{1}\join y_{2}) \\
			 & =x_{1}\join y_{1} && \text{ as }x_{1}\join y_{2}\in S_{1}\cap 
			 S_{2}. 
		\end{align*}
		Thus $x_{1}\geq y_{1}$. Dually $x_{1}\le y_{1}$. 
		
		In a similar way we have $x_{2}=y_{2}$. 
	
		\item 
		\begin{enumerate}[(a)]
			\item This part we do by localizing. Let $x, y\in\mathcal M$ and 
			let $\mathcal L=\mathcal M_{x\caret y}$.  Then we have 
			$\rsf C(\mathcal L)= [[x]\meet[y], \one]$ so we let 
			$\rsf G_{\mathcal L}=\rsf G\cap[[x]\meet[y], \one]= [[g], \one]$ 
			for some $g\geq x\caret y$. Also we note that 
			$S_{1, \mathcal L}= \eta_{\mathcal M}^{-1}[\rsf G_{\mathcal 
			L}]= S_{1}\cap\mathcal L = \mathcal L_{g}$, and 			
			$S_{2}\cap\mathcal L=\mathcal L_{g\to(x\caret y)}$. Let 
			$h=g\to(x\caret y)$ and $a=x\caret y$. 
			
			Thus we have for any $z\in\mathcal L$ that 
			$$
				z=(z\join\Delta(z\join g, g))\meet(z\join\Delta(z\join h, h))
			$$
			is the decomposition given by part (2) above. 
			Let $\phi_{g}=\phi_{\rsf G}\restrict\mathcal L$ so that 
			\begin{align*}
				\phi_{\rsf G}(w)&=(w\join\Delta(x\join g, g))\meet\Delta(\one, 
				w\join\Delta(x\join h, h))
			\end{align*}
			for any $w\in\mathcal L$. 
			
			We know that in any cubic implication algebra, if $a\sim b$ then 
			$f_{ab}\colon[a, \one]\to[b, \one]$ defined by 
			$$
				f_{ab}(w)=(w\join b)\meet(\Delta(\one, w)\join b)
			$$
			extends to an inner automorphism of $\mathcal L_{a}$ by 
			$$
				\widehat f_{ab}(z)=f_{ab}(z\join a)\meet \Delta\Bigl(\one, 
				f_{ab}\bigl((\Delta(\one, z)\join 
				a)\to a\bigr)\to b\Bigr). 
			$$
			
			Let 
			$$
				b=\Delta(g, a). 
			$$
			Then we have 
			\begin{tabular}[t]{>{$}r<{$}!{$=$}>{$}l<{$}} 
				b & \Delta(g, a)  \\
				 & g\meet\Delta(\one, g\to a)  \\
				 & (g\join a)\meet\Delta(\one, h\join a)=\phi_{\rsf G}(a)
			\end{tabular}\\
			is the image of the mapping restricted to $\mathcal L$. 
			
			We claim that $\phi_{g}$ is exactly the  inner 
			automorphism $\widehat f_{ab}$. 
			By \lemref{lem:onto} 
			it suffices to show that
			\begin{align}
				\widehat f_{ab}(z)\join\Delta(g, a) & =\phi_{g}(z)\join\Delta(g, a)
				\label{eq:oneAA}  \\
				\widehat f_{ab}(z)\join\Delta(\one, \Delta(g, a)) & =\phi_{g}(z)\join\Delta(\one, \Delta(g, a))
				\label{eq:twoAA}
			\end{align}
			for all $z$. Since we are working in an interval algebra, we will 
			use intervals from a Boolean algebra to do this. 
			
			Let $z=[z_{0}, z_{1}]$. Wolog $a=[0, 0]$, $g=[0, g]$ and $h=[0, 
			\comp g]$. Then $b=\Delta(g, a)=[g, g]$. As above 
			(\lemref{lem:intComp}) we have 
			\begin{align*}
				z\join\Delta(g\join z, g) & =[\comp g\meet z_{0}, g\join z_{1}]  \\
				z\join\Delta(h\join z, h) & =[g\meet z_{0},\comp  g\join z_{1}] \\
				\intertext{and so }
				\phi_{g}(z)&=[\comp g\meet z_{0}, g\join z_{1}]\meet[g\meet\comp{z}_{1}, 
				\comp g\join \comp z_{0}]\\
				&=[(\comp g\meet z_{0})\join(g\meet\comp z_{1}), (g\join 
				z_{1})\meet(\comp g\join \comp z_{0})]\\
				\phi_{g}(z)\join\Delta(g, a)&=\phi_{g}(z)\join[g, g]=[g\meet\comp z_{1}, g\join z_{1}]\\
				\phi_{g}(z)\join\Delta(\one, \Delta(g, a))&=\phi_{g}(z)\join[\comp g,\comp g]=[\comp g\meet z_{0}, \comp g\join \comp z_{0}]. 
			\end{align*}
			By general theory we have 
			\begin{align*}
				\widehat f_{ab}(z)\join\Delta(g, a) & =f_{ab}(z\join a)\\
				\widehat f_{ab}(z)\join\Delta(\one, \Delta(g, a)) & = \Delta\Bigl(\one, 
				f_{ab}\bigl((\Delta(\one, z)\join 
				a)\to a\bigr)\to b\Bigr). 
			\end{align*}
			This gives us 
			\begin{align*}
				\widehat f_{ab}(z)\join\Delta(g, a) & =f_{ab}(z\join a)\\
				 & =(z\join a\join b)\meet(\Delta(\one, z)\join\Delta(\one, 
				 a)\join b)  \\
				 & =([z_{0}, z_{1}]\join[0, g])\meet([\comp z_{1}, 
				 \comp z_{0}]\join[1, 1]\join[g, g]) \\
				 & =[0, z_{1}\join g]\meet[\comp z_{1}\meet g, \one]  \\
				 & =[\comp z_{1}\meet g, z_{1}\join g] \\
				 & =\phi_{g}(z)\join\Delta(g, a).  \\
			f_{ab}((\Delta(\one, z)\join 
				a)\to a) & = f_{ab}([0, \comp z_{0}]\to[0, 0])\\
				&=f_{ab}([0, z_{0}])\\
				&=([0, z_{0}]\join[g, g])\meet([\comp z_{0}, 1]\join[g, g])\\
				&=[0, z_{0}\join g]\meet[\comp z_{0}\meet g, 1]\\
				&=[\comp z_{0}\meet g, z_{0}\join g]\\
				\widehat f_{ab}(z)\join\Delta(\one, \Delta(g, a)) & = \Delta(\one, 
				f((\Delta(\one, z)\join 
				a)\to a)\to b)\\
				&=\Delta(\one, [\comp z_{0}\meet g, z_{0}\join g]\to[g, g])\\
				&=\Delta(\one, [z_{0}\meet g, \comp z_{0}\join g])\\
				&=[z_{0}\meet\comp g, \comp z_{0}\join\comp g]\\
				&=\phi_{g}(z)\join\Delta(\one, \Delta(g, a)). 
			\end{align*}
			
			It follows from this that $\phi_{\rsf G}$ is a well-defined, one-one, 
			and onto cubic homomorphism. 
			\begin{description}
				\item[Well-defined] Immediate from 
				the definition and part (2). 
			
				\item[One-one] As if $\phi_{\rsf G}(x)=\phi_{\rsf G}(y)$ we can 
				work as above and see that $\widehat f_{ab}(x)=\widehat 
				f_{ab}(y)$ so that $x=y$. 
			
				\item[{Homomorphism}] Let $x, y$ be given, and use $\mathcal L$ be 
				constructed with $x, y$. Then $x\join y$ and $\Delta(x, y)$ are in 
				$\mathcal L$ so we have 
				$\phi_{\rsf G}(x\join y)= \widehat f_{ab}(x\join y)= \widehat 
				f_{ab}(x)\join\widehat f_{ab}(y)= \phi_{\rsf G}(x)\join\phi_{\rsf 
				G}(y)$. Likewise $\Delta$ is preserved. 
				
				\item[Onto] To show this we note that $\phi_{\rsf 
				G}^{2}=\text{id}$ as 
				if $x=x_{1}\meet x_{2}$ with $x_{i}\in S_{i}$ then 
				the unique representation of $\phi_{\rsf G}(x)$ is
				$x_{1}\meet\Delta(\one, x_{2})$. Thus
				\begin{align*}
					\phi_{\rsf G}^{2}(x) & =\phi_{\rsf G}(x_{1}\meet\Delta(\one, x_{2}))  \\
					 &  = x_{1}\meet\Delta(\one, \Delta(\one, x_{2})) \\
					 & =x_{1}\meet x_{2}=x. 
				\end{align*}
				Since this is true, if $x\in\mathcal M$ then $x=\phi_{\rsf 
				G}^{2}(x)$ is in the range of $\phi_{\rsf G}$. 
			\end{description}
		
			\item We can argue as above to see that for any 
			$x\in\mathcal M$ we have 
			$\phi_{\rsf G}(x)= \widehat f_{ab}(x)\sim x$ so that $\rsf C(\phi_{\rsf 
			G})=\text{id}$. 
		
			\item If $x\in S_{1}$ then $x=x\meet \one$ so that 
			$\phi_{\rsf G}(x)= x\meet\Delta(\one, \one)= x$. 
			
			If $\phi_{\rsf G}(x)= x= x_{1}\meet x_{2}$ with $x_{i}\in S_{i}$ then we have 
			$x_{1}\meet x_{2}= x_{1}\meet \Delta(\one, x_{2})$. By the 
			uniqueness of representation we have $x_{2}= \Delta(\one, x_{2})$ 
			and so $x_{2}=\one$. Thus $x=x_{1}\in S_{1}$. 
		\end{enumerate}
	\end{enumerate}
\end{proof}

\begin{thm}\label{thm:isoGroups}
	$$
		\inn(\mathcal M)\simeq\inn(\rsf I\rsf C(\mathcal M)). 
	$$
\end{thm}
\begin{proof}
	We know from \cite{BO:fil} that $\inn(\rsf I\rsf C(\mathcal M))$ is 
	isomorphic to the group of $\rsf C(\mathcal M)$-Boolean filters. 
	
	We define a mapping $\Omega$ from $\inn(\mathcal M)$ to this group by
	$$
		\Omega(\phi)=\rsf C(\mathcal M_{\phi}). 
	$$
	By \lemref{lem:BoolCC} this filter is $\rsf C(\mathcal M)$-Boolean. 
	By \lemref{lem:collapseDewt} and \thmref{thm:MPhiIsGood} this 
	mapping is one-one. By the last theorem the mapping is onto. We just 
	need to show that it is a homomorphism, i.e. that
	$$
		\rsf C(\mathcal M_{\phi_{1}\phi_{2}})=\rsf C(\mathcal M_{\phi_{1}})+\rsf C(\mathcal M_{\phi_{2}}). 
	$$
	First notice that $\mathcal M_{\phi_{1}\phi_{2}}=\Set{x | 
	\phi_{1}(x)=\phi_{2}(x)}$ -- since both maps are their own inverses. 
	By definition
	\begin{multline*}
		\rsf C(\mathcal M_{\phi_{1}})+\rsf C(\mathcal M_{\phi_{2}})=
		\Bigl[\bigl(\rsf C(\mathcal M_{\phi_{1}})\to\rsf C(\mathcal 
		M)\bigr)\cap\bigl(\rsf C(\mathcal 
		M_{\phi_{2}})\to\rsf C(\mathcal M)\bigr)\Bigr]\join\\
		\Bigr[\rsf C(\mathcal M_{\phi_{1}})\cap\rsf C(\mathcal 
		M_{\phi_{2}})\Bigl]. 
	\end{multline*}
	Let $x\in\mathcal M$ be arbitrary. 
	Then $x=x_{1}\meet x_{2}$ for some $x_{1}\in \mathcal M_{\phi_{1}}$ and 
	$x_{2}\in D_{\phi_{1}}$. Now let
	$x_{1}=x_{11}\meet x_{12}$ and $x_{2}=x_{21}\meet x_{22}$ where 
	$x_{i1}\in \mathcal M_{\phi_{2}}$ and 
	$x_{i2}\in D_{\phi_{2}}$. Then we have 
	\begin{align*}
		\phi_{1}(x) & =x_{11}\meet x_{12}\meet\Delta(\one, x_{21})\meet 
		\Delta(\one, x_{22})  \\
		 \phi_{2}(x) & =x_{11}\meet \Delta(\one, x_{12})\meet x_{21}\meet 
		\Delta(\one, x_{22}). 
	\end{align*}
	Then if $\phi_{1}(x)=\phi_{2}(x)$ the uniqueness of the 
	representations implies that $x_{12}=\Delta(\one, x_{12})$ and 
	$x_{21}=\Delta(\one, x_{21})$. Thus $x_{12}=x_{21}=\one$. 
	Therefore $x=x_{11}\meet x_{22}$ with $x_{11}\in\mathcal 
	M_{\phi_{1}}\cap\mathcal M_{\phi_{2}}$ and 
	$x_{22}\in D_{\phi_{1}}\cap D_{\phi_{2}}$. As this entails 
	$x=x_{11}\caret\Delta(\one, x_{22})$ we have 
	$[x]=[x_{11}]\meet[x_{22}]$ and $[x_{11}]\in \Bigr[\rsf C(\mathcal M_{\phi_{1}})\cap\rsf C(\mathcal 
			M_{\phi_{2}})\Bigl]$, 
			$[x_{22}]\in\Bigl[\bigl(\rsf C(\mathcal M_{\phi_{1}})\to\rsf C(\mathcal 
		M)\bigr)\cap\bigl(\rsf C(\mathcal 
		M_{\phi_{2}})\to\rsf C(\mathcal M)\bigr)\Bigr]$. Thus
		$[x]\in \rsf C(\mathcal M_{\phi_{1}})+\rsf C(\mathcal M_{\phi_{2}})$. 
		
	Conversely, if $[x]\in \rsf C(\mathcal M_{\phi_{1}})+\rsf C(\mathcal 
	M_{\phi_{2}})$ implies $[x]=[x_{1}]\meet[x_{2}]$ for some 
	$[x_{1}]$ in $\Bigr[\rsf C(\mathcal M_{\phi_{1}})\cap\rsf C(\mathcal 
			M_{\phi_{2}})\Bigl]$, 
			$[x_{2}]$ in $\Bigl[\bigl(\rsf C(\mathcal M_{\phi_{1}})\to\rsf C(\mathcal 
		M)\bigr)\cap\bigl(\rsf C(\mathcal 
		M_{\phi_{2}})\to\rsf C(\mathcal M)\bigr)\Bigr]$, 
		and (as in previous analyses) we may assume that 
		$x=x_{1}\meet x_{2}$ where $x_{1}\in \mathcal 
	M_{\phi_{1}}\cap\mathcal M_{\phi_{2}}$ and $x_{2}\in D_{\phi_{1}}\cap D_{\phi_{2}}$.
	Then we have 
	\begin{align*}
		\phi_{1}(x) & =x_{1}\meet \Delta(\one, x_{2})  \\
		\intertext{and }
		 \phi_{2}(x) & =x_{1}\meet \Delta(\one, x_{2}). 
	\end{align*}
	Thus $[x]\in\rsf C(\mathcal M_{\phi_{1}\phi_{2}})$. 
\end{proof}

The isomorphism we have constructed in this theorem comes about in a 
rather indirect fashion. This is in sharp contrast to earlier 
isomorphism results that came about via extension of homomorphism 
results. In this case such extensions seem impossible to obtain. 

% The title of this section included the words ``Strong Congruences'',  
% but so far we have avoided any cubic congruences. The cubic congruences come in 
% because we know from \cite{BO:Cong} that the lattice of cubic congruences is 
% isomorphic to the lattice of upwards-closed MR-subalgebras. The theorems 
% above have established a relationship between inner automorphisms and 
% upwards-closed subalgebras. The particular congruences that inner 
% automorphisms correspond to are \emph{strong congruences}. More on 
% them later. 

\end{document}